\documentclass{amsproc}
\usepackage{amsfonts}
\usepackage{amssymb, amsmath, 
}

\usepackage{epsf}            
\usepackage{amssymb, amsmath}
\usepackage {graphicx}
\usepackage {color}
\usepackage{hhline}

\begin{document}

\newcommand{\bq}{\begin{equation}}
\newcommand{\eq}{\end{equation}}
\newcommand{\rf}[1]{(\ref{#1})}

\def\bkN{{\rm I\kern-.20em N}}
\def\bkR{{\rm I\kern-.17em R}}
\def\bkK{{\rm I\kern-.22em K}}
\def\bkC{{\rm \kern.24em
       \vrule width.05em height1.4ex depth-.05ex
       \kern-.26em C}}
\def\BBox{\vrule height 0.5em width 0.6em depth 0em}

\numberwithin{equation}{section}

\newcommand{\ZZ}{{\mathbb Z} }
\newcommand{\NN}{{\mathbb N}}
\newcommand{\CC}{{\mathbb C}}
\newcommand{\btd}{\nabla}
\newcommand{\btu}{\Delta}
\newcommand{\dst}{\displaystyle}
\newcommand{\half}{{1/2}}
\newcommand{\esp}{\quad}
\newcommand{\fhyp}{\textrm{\large F}}
\newtheorem{teo}{Theorem}[section]
\newtheorem{cor}[teo]{Corollary}
\newtheorem{lema}[teo]{Lemma}
\newtheorem{prop}[teo]{Proposition}
\title{Basic Fourier series: convergence on and outside the $q$-linear grid}
\author{J. L. Cardoso}



\begin{abstract}
A \emph{$\,q$-type H\"older condition} on a function
$\,f\,$ is given in order to establish (uniform) convergence of the
corresponding basic Fourier series $\,S_q[f]\,$ to the function
itself, on the set of points of the $\,q$-linear grid.

Furthermore, by adding others conditions, one guaranties the (uniform)
convergence of $\,S_q[f]\,$ to $\,f\,$ on and "outside" the set points of
the $\,q$-linear grid.
\end{abstract}

\maketitle


AMS MOS: 33D15, 42C10


\section{\sc{Introduction}}

Basic Fourier expansions on $\,q$-quadratic and on $\,q$-linear
grids were first considered in \cite{BS:1998} and in \cite{BC:2001},
respectively. Recently, in \cite{C:2005}, sufficient conditions for (uniform)
convergence of the $\,q$-Fourier series in terms of basic trigonometric
functions $\,S_q\,$ and $\,C_q\,$, on a $\,q$-linear grid, were given.
In \cite{S:1997} it was established an "addition" theorem
for the corresponding basic exponential function, being these functions
equivalent to the ones introduced by H. Exton in \cite{E:1983}.
Following the unified approach of M. Rahman in \cite{R:1998}, these
functions can be seen as analytic linearly independent solutions of
the initial value problem
\begin{equation*}
\frac{\delta f(x)}{\delta x}=\lambda f(x)\,,\quad f(0)=1\,,
\end{equation*}
where $\,\delta\,$ is the symmetric $\,q$-difference operator acting on
a function $\,f\,$ by
\bq \delta f(x)=f(q^{1/2}x)-f(q^{-1/2}x)\,, \label{1.1} \eq
with $\,0<q<1\,.$
Then, from \rf{1.1}, \bq \frac{\delta f(x)}{\delta
x}=\frac{f(q^{1/2}x)-f(q^{-1/2}x)}{x(q^{1/2}-q^{-1/2})}\:.
\label{1.2} \eq
There exists an important relation between this
difference operator and the $\,q$-integral. The $\,q$-integral is
defined by
$$\int_{0}^{a}f(x)d_qx=a(1-q)\sum_{n=0}^{\infty}f(aq^n)q^n$$
and
\begin{equation}
\int_{a}^{b}f(x)d_qx=\int_{0}^{b}f(x)d_qx-\int_{0}^{b}f(x)d_qx
\label{1.3}\,.
\end{equation}
From \rf{1.2} and \rf{1.3} it follows
\begin{equation}
\int_{-1}^1\frac{\delta f(x)}{\delta x}d_qx=q^{\frac 1 2}\left\{\left[
f(q^{-\frac 1 2})-f(-q^{-\frac 1 2})\right]-\left[f(0^+)-f(0^-)\right]\right\}\,,
\label{1.4}
\end{equation}
hence,
one have the following formula \cite{C:2005} for $\,q$-integration by parts:
\begin{equation}
\begin{array}{l}
\displaystyle\int_{-1}^1g\big(q^{\pm\frac 1 2}x\big)\frac{\delta_q f(x)}{\delta_q x}d_qx
\:=\:-\int_{-1}^1f\big(q^{\mp\frac 1 2}x\big)\frac{\delta_q g(x)}{\delta_q x}d_qx\:+
\\ [1em]
\displaystyle\hspace{1.1em}
q^{\frac 1 2}\left\{\left[\,\big(fg\big)\big(q^{-\frac 1 2}\big)-\big(fg\big)
\big(-q^{-\frac 1 2}\big)\,\right]-\left[\,\big(fg\big)\big(0^+\big)-
\big(fg\big)\big(0^-\big)\,\right]\right\}.
\end{array}
\label{1.5}
\end{equation}

These functions satisfy an orthogonality relation \cite{BC:2001, E:1983}
where the corresponding inner product is defined in terms of the $\,q$-integral \rf{1.4}. In
\cite{BC:2001}, it was proved that they form a complete system and
analytic bounds on their roots were derived.

As we will refer in section 2, the above $\,q$-trigonometric
functions can be written using the Third Jackson $\,q$-Bessel
funtion (or the Hahn-Exton $\,q$-Bessel function). In
\cite{ABC:2003}, analytic bounds were derived for the zeros of
this function --which includes, as particular cases, the
corresponding results established in \cite{BC:2001}-- and
recently, in \cite{AB:2005}, it was shown that they define a
complete system.

Throughout this paper we will follow the notation used in \cite{GR:1997}
which is now standard.

The publications \cite{BC:2001, BS:1998, C:2000, C:2005, 
S:2002, S:2003}
are the most affiliated with this work.
For other type of expansions (sampling theory) or related topics see
\cite{A:2005, A:2005b, A:2006, ABC:2003, AM:2005}.

\section{\sc The $\,q$-Linear Sine and Cosine. Properties.}

The initial value problem
\begin{equation*} 
\frac{\delta f(x)}{\delta x}=\lambda
f(x)\:,\quad f(0)=1\,,
\end{equation*}
has the analytic solution \cite{BC:2001}
\bq
\exp_{q}[\lambda(1-q)z]=\sum_{n=0}^{\infty}\frac{[\lambda(1-q)z]^nq^{(n^2-n)/4}}{(q;q)_n}\:,
\label{2.2}
\eq
which is a standard $q$-analog of the classical exponential function
\cite{GR:1997,R:1998}. The $\,q$-linear sine and cosine, $\:S_q(z)\:$ and
$\:C_q(z)\,$, are then defined by $$\exp_{q}{iz}:=C_q(z)+iS_q(z)\:.$$
From \rf{2.2} we get
\begin{equation*} 
C_q(z)=\sum_{n=0}^{\infty}\frac{(-1)^nq^{n[n-(1/2)]}z^{2n}}{(q;q^2;q^2)_n}\:=
\:{_1\phi_1}\left(\begin{array}{ccc}
 \! 0 \! & \!\! & \!\! \\ [-0.5em]
 \!\! & \! ;\! & \! q^2,\:q^{1/2}z^2 \! \\ [-0.5em]
 \! q \! & \!\! & \!\! \,,
\end{array}\right)
\end{equation*}
\begin{equation*} 
S_q(z)=\frac{z}{1-q}\sum_{n=0}^{\infty}\frac{(-1)^nq^{n[n+(1/2)]}z^{2n}}{(q^2;q^3;q^2)_n}\:=
\:\frac{z}{1-q}\:{_1\phi_1}\left(\begin{array}{ccc}
 \! 0 \! & \!\! & \!\! \\ [-0.5em]
 \!\! & \! ;\! & \! q^2,\:q^{3/2}z^2 \! \\ [-0.5em]
 \! q^3 \! & \!\! & \!\!
\end{array}\right)\,,
\end{equation*}
which can be written in terms of the third Jackson $\,q$-Bessel function
(or, Hahn-Exton $\,q$-Bessel function) \cite{I:1982, KS:1994, S:1992}
\begin{equation*} 
J_{\nu}^{(3)}(z;q):=z^{\nu}\frac{\left(q^{\nu+1};q\right)_{\infty}}{(q;q)_{\infty}}
\:{_1\phi_1}\left(\begin{array}{ccc}
 \! 0 \! & \!\! & \!\! \\ [-0.5em]
 \!\! & \! ;\! & \! q,\:qz^2 \! \\ [-0.5em]
 \! q^{\nu+1} \! & \!\! & \!\!
\end{array}\right)
\end{equation*}
as
\begin{equation*} 
C_q(z)=q^{-3/8}\frac{(q^2;q^2)_{\infty}}{(q;q^2)_{\infty}}
z^{1/2}J_{-1/2}^{(3)}\left(q^{-3/4}z;q^2\right)\,,
\end{equation*}
\begin{equation*} 
S_q(z)=q^{1/8}\frac{(q^2;q^2)_{\infty}}{(q;q^2)_{\infty}}
z^{1/2}J_{1/2}^{(3)}\left(q^{-1/4}z;q^2\right)
\end{equation*}
They satisfy \cite{BC:2001}
\bq \label{2.8}
\frac{\delta C_q(\omega z)}{\delta
z}=-\frac{\omega}{1-q}S_q(\omega z)\,,
\eq
\bq \label{2.9}
\frac{\delta S_q(\omega z)}{\delta z}=\frac{\omega}{1-q}C_q(\omega z)\,,
\eq
and, when $\,\omega\,$ is such that $\,S_q(\omega)=0\,$,
\bq \label{2.10}
\big[C_q(\omega)\big]^{-1}=C_q(q^{-1/2}\omega)=C_q(q^{1/2}\omega)\,.
\eq
It is known \cite{BC:2001} that the roots of $\,C_q(z)\,$ and
$\,S_q(z)\,$ are {\em real}, {\em simple} and {\em countable}.
Further, because $\,C_q(z)\,$ and $\,S_q(z)\,$ are
respectively even and odd functions, the roots of
$\,C_q(z)\,$ and $\,S_q(z)\,$ are symmetric and we will denote the
positive zeros of the function $\,S_q(z)\,$ by $\,\omega_k\,,$
$k=1,2,\ldots$, with $\,\omega_1<\omega_2<\omega_3<\ldots$.

As we mentioned before, the zeros of the function $\:S_q(z)\:$
form a discrete set of symmetric points in the real line.
In \cite[page 145]{BC:2001}, it was shown that the set of positive
zeros $\:\omega_k\,,$ $k=1,2,\ldots$ of the function $\,S_q(z)\,,$
verify the following {\em analytic bounds}:

\vspace{0.3em}
\textsl{
If $\:0<q<\beta_0\,,$ where $\:\beta_0\:$ is the root of
$\:(1-q^2)^2-q^3\,,$ $\:0<q<1\,,$ then
$$q^{-k+\alpha_k+1/4}<\omega_k<q^{-k+1/4}\:,\quad k=1,2,\ldots\:,$$
\hspace{1.6em}where
$$\alpha_k\equiv\alpha_k(q)=\frac{\log{\left[1-\frac{q^{2k+1}}{1-q^{2k}}\right]}}{2\log{q}}\:,
\quad k=1,2,\ldots\:.\:$$
}
According to {\em Remark 1} in \cite[page 145]{BC:2001}, the previous
result can be restated in the following form:

\vspace{0.5em}
\noindent\textbf{Theorem A}\textsl{
For every $\:q\:$, $\:0<q<1\,,$ $\:K\:$ exists such that if $\:k\geq K\:$ then
\begin{equation*}
\omega_k=q^{-k+\epsilon_k+1/4}\:,\quad 0<\epsilon_k<\alpha_k(q)\:.
\end{equation*}
}

By using Taylor expansion one finds out that
\begin{equation}\label{2.11}
\alpha_k(q)=\mathcal{O}(q^{2k})\quad\mbox{as}\quad k\rightarrow\infty\,.
\end{equation}

\vspace{0.7em}
Theorem 4.1 of \cite[page 139]{BC:2001} settle the {\em orthogonality relations}:

\vspace{0.7em}
\noindent\textbf{Theorem B}\textsl{
Considering
$\,\mu_k=(1-q)C_q(q^{1/2}\omega_k)S_q^{\prime}(\omega_k)\,$
we have
\begin{equation*}
\int_{-1}^1C_q(q^{1/2}\omega_k x)C_q(q^{1/2}\omega_m x)d_qx=
\left\{\begin{array}{lcl} 0 & \mbox{if} &
k\neq m
\\ 2 & \mbox{if} & k=0=m \\
\mu_k & \mbox{if} & k=m\neq 0
\end{array}\right.
\end{equation*}
\begin{equation*}
\int_{-1}^1S_q(q\omega_k x)S_q(q\omega_m x)d_qx=
\left\{\begin{array}{lcl} 0 & if &
k\neq m\vee k=0=m \\
 & & \\
q^{-1/2}\mu_k & if & k=m\neq 0
\end{array}\right.\,.
\end{equation*}
}

\vspace{0.7em}
The \emph{Completeness Theorem} \cite[page 153]{BC:2001}, where a misprint is corrected,
states the following:

\vspace{0.5em}
\noindent\textbf{Theorem C}\textsl{
Let $\:f(\omega_k z)=C_q\!\left(q^{\frac{1}{2}}\omega_k z\!\right)+
iS_q\!\left(q\omega_k z\right)\:$ where the $\,\omega_k\,,$ $\,\omega_0\!=\!0<\!\omega_1\!
<\!\omega_2\!<\!\ldots\,$ are the non-negative roots of $\,S_q(z)\,$.
Suppose that
$$\int_{-1}^{1}g(z)f(\omega_kz)d_qz=0\quad,\hspace{2em} k=0,1,2,\ldots$$
where $\,g(z)\,$ is bounded on
$\:z=\pm q^j\:,\;j=0,1,2,\ldots\:.$
Then, $\:g(z)\equiv 0\:,$ i.e., $\:g\left(\pm q^j\right)=0\:$ for all
$\;j=0,1,2,\ldots\:.$
}

\vspace{0.7em}
To end this section we write down the Theorem 6.2 of \cite[page 150]{BC:2001}:

\vspace{0.5em}
\noindent\textbf{Theorem D}\textsl{
If $\;S_q(\omega_k)=0\;$ then, for
$\;n=0,1,2,\ldots\:,$
\begin{equation*}
S_q(q^{1+n}\omega_k)=S_q(q\omega_k)\sum_{j=0}^{n}(-1)^{j}q^{j(j+\frac{1}{2})}
\frac{\left(q^{1+n-j};q\right)_{2j+1}}{(q;q)_{2j+1}}\left(\omega_k^2\right)^j\,,
\end{equation*}
\begin{equation*}
C_q(q^{\frac 1 2+n}\omega_k)=C_q(q^{\frac 1 2}\omega_k)\sum_{j=0}^{n}
(-1)^{j}q^{j(j-\frac{1}{2})}\frac{\left(q^{1+n-j};q\right)_{2j}}{(q;q)_{2j}}
\left(\omega_k^2\right)^j\,.
\end{equation*}
}

\section{\sc The Fourier Coefficients}

As a consequence of the orthogonality relations of Theorem B, we may consider
formal Fourier expansions of the form
\bq
f(x)\sim S_q[f](x)=\frac{a_0}{2}+\sum_{k=1}^{\infty}\left[a_kC_q\left(q^{\frac{1}{2}}\omega_k
x\right)+ b_kS_q\left(q\omega_k x\right)\right] \:, \label{3.1}
\end{equation}
with $\:a_0=\int_{-1}^{1}f(t)d_qt\:$ and, for $\:k=1,2,3,\ldots\:,$
\begin{equation}
a_k=\frac{1}{\mu_k}\int_{-1}^{1}f(t)
C_q\left(q^{\frac{1}{2}}\omega_k t\right)d_qt \label{3.2} \eq \bq
b_k=\frac{q^{\frac{1}{2}}}{\mu_k}\int_{-1}^{1}f(t)
S_q\left(q\omega_k t\right)d_qt \:, \label{3.3} \eq
where
\bq
\mu_k=(1-q)C_q(q^{1/2}\omega_k)S_q^{\prime}(\omega_k)\,.\,
\label{3.4}
\eq
In order to study the convergence of the series (\ref{3.1})-(\ref{3.4}),
it becomes clear that we need to know the behavior of the factor $\,\mu_k\,$
of the denominator as $\:k\rightarrow\infty\,,$ which is equivalent to
control the behavior of $\:S_q^{\prime}(\omega_k)\:$ and
$\,C_q(q^{1/2}\omega_k)\,$ as $\:k\rightarrow\infty\,.$

\vspace{0.8em}
Theorem 3.2 from \cite{C:2005} asserts that

\vspace{0.6em}
\noindent\textbf{Theorem E}\textsl{\: At least for $\:0<q\leq (1/51)^{1/50}\,,$
$$S_q^{\prime}(\omega_k)=\frac{2}{1-q}q^{-(k-\frac 1
2-\epsilon_k)^2}S_k \:,$$ where $\:S_k\,$ satisfies
$\:\displaystyle\liminf_{k\rightarrow\infty}|S_k|>0\,.$
}

\vspace{1em}
With respect to $\,S_k\,$ from the previous theorem we have the following
lemma:
\begin{lema}\label{L3.1}
There exists a constant $\,B\,,$ independent of $\,k\,,$ such that
$$\left|S_k\right|\leq B\;,\quad k=1,2,3,\ldots\,.$$
\end{lema}
\begin{proof} The expression of $\,S_k\,$ is given \cite[page 147]{BC:2001} by
$$S_k=\sum_{n=0}^{\infty}\frac{(-1)^n
n q^{(n-k+1/2+\varepsilon_k)^2}}{(q^2,q^3;q^2)_n}=
(-1)^k\sum_{m=-k}^{\infty}\frac{(-1)^m
m q^{(m+1/2+\varepsilon_k)^2}}{(q^2,q^3;q^2)_{m+k}}\,.$$
For $\,k\,$ large enough, by Theorem A and (\ref{2.11}),
$\,1/2+\varepsilon_k>0\,$ hence
$$\left|S_k\right|\leq\sum_{m=-k}^{\infty}\frac{
|m| q^{(m+1/2+\varepsilon_k)^2}}{(q^2,q^3;q^2)_{m+k}}\leq
\frac{2}{(q^2;q)_{\infty}}\sum_{m=1}^{\infty}
m q^{(m-1)^2}=B$$
which completes the proof since the infinite series on the right member
is convergent.
\end{proof}

\noindent We observe that the constant $B$, as well as $S_k$, depend on
the parameter $q\,.$

\vspace{0.8em}
The behavior of $\,C_q(q^{1/2}\omega_k)\,$ as $\:k\rightarrow\infty\,$ will
be known by the corresponding behavior of $\:C_q(\omega_k)\:$ and by
(\ref{2.10}). Theorem 3.3 of \cite{C:2005} establishes

\vspace{0.6em}
\noindent\textbf{Theorem F}\textsl{\: At least for $\;\;0<q\leq (1/50)^{1/49}\,,$
$$C_q(\omega_k)=q^{-(k-\epsilon_k)^2}R_k \:,$$ \hspace{2em} where
$\quad\left|R_k\right|<\displaystyle\frac{2}{(1-q)(q;q)_{\infty}}\quad$
and
$\quad\displaystyle\liminf_{k\rightarrow\infty} |R_k|>0\,.$
}

\vspace{0.8em}
To end this section, we collect the Theorems 4.1, 4.2 and 4.3 of
\cite{C:2005}:

\vspace{0.6em}
\noindent\textbf{Theorem G}\textsl{
If $\:c\in\mathbb{R}\:$ exists such that, as $\:k\rightarrow\infty\,,$
$$\int_{-1}^1f(t)C_q\left(q^{\frac 1 2}\omega_kt\right)d_qt=\mathcal{O}
\left(q^{ck}\right)\quad\mbox{and}\quad\int_{-1}^1f(t)S_q\left(q\omega_kt\right)d_qt=
\mathcal{O}\left(q^{ck}\right)$$
then, at least for $\,0<q\leq (1/51)^{1/50}\,,$ the $\,q$-Fourier series
\rf{3.1} is pointwise convergent at each fixed point
$\,\displaystyle x\in V_q=\left\{\pm q^{n-1}:\,n\in\bkN\,\right\}\,.$
}

\vspace{0.8em}
\noindent\textbf{Theorem H}\textsl{
If $\:c>1\:$ exists such that, as $\:k\rightarrow\infty\,,$ \vspace{-0.4em}
$$\int_{-1}^1f(t)C_q\left(q^{\frac{1}{2}}\omega_kt\right)d_qt=\mathcal{O}
\left(q^{ck}\right)\quad\mbox{and}\quad\int_{-1}^1f(t)S_q\left(q\omega_kt\right)d_qt=
\mathcal{O}\left(q^{ck}\right)$$
\noindent then, the $\:q$-Fourier series
\rf{3.1}, at least for $\,0<q\leq (1/51)^{1/50}\,,$ converges uniformly
on $\,\displaystyle V_q=\left\{\pm q^{n-1}
:\,n\in\bkN\,\right\}\,.$
}

\vspace{0.8em}
\vspace{0.8em}
\noindent\textbf{Theorem I}\textsl{
If $\;f\,$ is a bounded function on the set
$\,\displaystyle \,V_q=\left\{\pm q^{n-1}
:\,n\in\bkN\,\right\},$ and the $\,q$-Fourier series $\,S_q[f](x)\,$
converges uniformly on $\,V_q\,$ then its sum is $\,f(x)\,$ whenever
$\,x\in V_q\,.$
}

\section{Convergence condition on the function}

Denoting the $\,q$-Fourier coefficients of a function $\,f\,$ by
$\,a_k\big(f(x)\big)\,$ and $\,b_k\big(f(x)\big)\,,$ $\,k=1,2,3,\ldots\,,$
using (\ref{3.2})-(\ref{3.4}) and (\ref{2.8})-(\ref{2.9}) one have, by (\ref{1.5}),
\begin{equation}\label{4.1}
a_k\big(f(x)\big)-\frac{1-q}{q^{1/2}\omega_k\mu_k}\int_{-1}^{1}S_q\left(q\omega_k t\right)
\frac{\delta f\big(q^{\frac 1 2}t\big)}{\delta t}d_qt-\frac{1-q}{q\omega_k}\:
b_k\!\left(\frac{\delta f(q^{\frac 1 2}x)}{\delta x}\right)
\end{equation}
and
\begin{equation}\label{4.2}
\begin{array}{l}
b_k(f(x))=\displaystyle\frac{q-1}{q^{\frac 1 2}\omega_k\mu_k}\:\Big\{
q^{\frac 1 2}\Big[f\big(q^{-1}\big)-f\big(-q^{-1}\big)\Big]
C_q\left(q^{\frac 1 2}\omega_k\right)- \\
\displaystyle\hspace{8em}q^{\frac 1 2}\Big[f\big(0^+\big)-f\big(0^-\big)\Big]-
\int_{-1}^{1}C_q\left(q^{\frac 1 2}\omega_k t\right)
\frac{\delta f\big(q^{-\frac 1 2}t\big)}{\delta t}d_qt\:\Big\} \\ [1.5em]
=\displaystyle
\frac{1-q}{q^{\frac 1 2}\omega_k}\,\left\{
a_k\!\left(\frac{\delta f(q^{-\frac 1 2}x)}{\delta x}\right)+
q^{\frac 1 2}\left[\frac{f\big(0^+\big)\!-\!f\big(0^-\big)}{\mu_k}-
\frac{f\big(q^{-1}\big)\!-\!f\big(\!-q^{-1}\big)}{(1-q)S_q^{\prime}(\omega_k)}\right]
\right\}.
\end{array}
\end{equation}
The conjugation of this last two identities with Theorem H enables us to deduce
conditions on the function $\,f\,$ in order to guarantee uniform convergence of
the corresponding Fourier series $\,S_q[f]\,$.
In its statement, we will consider the notation

\vspace{-0.6em}
$$L_q^{\infty}[-1,1]=\Big\{f:\sup\big\{\left|f\big(\pm q^{n-1}\big)\right|:
n\in\mathbb{N}\big\}<\infty\Big\}$$
and the following definition:

\vspace{0.7em}
\noindent$\mathbf{Definition\:4.1}$ \textsl{If two constants $\,M\,$ and $\,\lambda\,$ exist such that
\begin{equation}\label{4.3}
\Big|f\big(\pm q^{n-1}\big)-f\big(\pm q^{n}\big)\Big|\leq M q^{\lambda n}\:,\quad
n=0,1,2,\ldots\,,
\end{equation}
then the function $\,f\,$ is said to be $\,q$-linear H\"older of order
$\,\lambda\,.$}

\begin{teo}\label{T4.1}
If $\,f\in L_q^{\infty}[-1,1]\,$ is a $\,q$-linear H\"older function of order
$\,\lambda>\frac 1 2\,$ and satisfies $\,f(0^+)=f(0^-)\,$
then, at least for $\,0<q\leq (1/50)^{1/49}\,,$ the corresponding $\,q$-Fourier
series $\,S_q[f]\,$ converges uniformly to $\,f\,$ on the set of points
$\,\displaystyle V_q=\left\{\pm q^{n-1}:\,n\in\bkN\,\right\}\,.$
\end{teo}
\begin{proof} From (\ref{3.2}) and (\ref{4.1}) one have
\begin{equation}\label{4.4}
\int_{-1}^1f(t)C_q\left(q^{\frac 1 2}\omega_kt\right)d_qt=\mu_k\,a_k\big(f\big)
=-\frac{1-q}{q^{1/2}\omega_k}\int_{-1}^{1}S_q\left(q\omega_k t\right)
\frac{\delta f\big(q^{\frac 1 2}t\big)}{\delta t}d_qt\,.
\end{equation}
Similarly, from (\ref{3.3}) and (\ref{4.2}),
\begin{equation}\label{4.5}
\begin{array}{l}
\displaystyle\int_{-1}^1f(t)S_q\left(q\omega_kt\right)d_qt =
\displaystyle q^{-1/2}\mu_k\,b_k\big(f\big)= \\ [1em]
\displaystyle
\frac{q\!-\!1}{q\omega_k}\:\Big\{
q^{\frac 1 2}\Big[f\big(q^{-1}\big)\!-\!f\big(-q^{-1}\big)\Big]
C_q\left(q^{\frac 1 2}\omega_k\!\right)- \int_{-1}^{1}\!C_q\!\left(q^{\frac 1 2}\omega_k t\right)
\frac{\delta f\big(q^{-\frac 1 2}t\big)}{\delta t}d_qt\Big\}. \\ [1em]
\end{array}
\end{equation}
By Cauchy-Schwarz inequality we have
\begin{equation}\label{4.6}
\!\left|\int_{-1}^{1}\!S_q\left(q\omega_k t\right)
\frac{\delta f\big(q^{\frac 1 2}t\big)}{\delta t}d_qt\right|\leq
\left(\int_{-1}^{1}\!
S_q^2\left(q\omega_k t\right)d_qt\right)^{\frac 1 2}
\left(\int_{-1}^{1}\!\left(\frac{\delta f\big(q^{\frac 1 2}t\big)}{\delta t}
\right)^2d_qt\right)^{\frac 1 2}
\end{equation}
and
\begin{equation}\label{4.7}
\left|\int_{-1}^{1}\!\!C_q\!\left(q^{\frac 1 2}\omega_k t\right)
\frac{\delta f\big(q^{-\frac 1 2}t\big)}{\delta t}d_qt\right|\leq
\displaystyle\left(\int_{-1}^{1}\!\!
C_q^2\!\left(q^{\frac 1 2}\omega_k t\right)d_qt\right)^{\!\frac 1 2}
\!\left(\int_{-1}^{1}\!\!\left(\frac{\delta f\big(q^{-\frac 1 2}t\big)}{\delta t}
\right)^2\!\!d_qt\!\right)^{\!\frac 1 2}
\end{equation}
Using the orthogonality relations of Theorem B we may write
\begin{equation*} 
q^{\frac 1 2}\int_{-1}^{1}S_q^2\left(q\omega_k t\right)d_qt=
\int_{-1}^{1}C_q^2\left(q^{\frac 1 2}\omega_k t\right)d_qt=
\mu_k=(1-q)C_q\left(q^{\frac 1 2}\omega_k\right)S_q^{\prime}(\omega_k)\,,
\end{equation*}
thus (\ref{4.6}) and (\ref{4.7}) become, respectively,
\begin{equation}\label{4.9}
\begin{array}{l}
\displaystyle\left|\int_{-1}^{1}S_q\left(q\omega_k t\right)
\frac{\delta f\big(q^{\frac 1 2}t\big)}{\delta t}d_qt\right|\leq \\ [1em]
\hspace{3.5em}\displaystyle q^{-\frac 1 4}(1-q)^{\frac 1 2}\left(C_q
\left(q^{\frac 1 2}\omega_k\right)S_q^{\prime}(\omega_k)\right)^{\frac 1 2}
\left(\int_{-1}^{1}\left(\frac{\delta f\big(q^{\frac 1 2}t\big)}{\delta t}
\right)^2d_qt\right)^{1/2}
\end{array}
\end{equation}
and
\begin{equation}\label{4.10}
\begin{array}{l}
\displaystyle\left|\int_{-1}^{1}C_q\left(q^{\frac 1 2}\omega_k t\right)
\frac{\delta f\big(q^{-\frac 1 2}t\big)}{\delta t}d_qt\right|\leq \\ [1em]
\hspace{3.5em}\displaystyle(1-q)^{\frac 1 2}\left(C_q
\left(q^{\frac 1 2}\omega_k\right)S_q^{\prime}(\omega_k)\right)^{\frac 1 2}
\left(\int_{-1}^{1}\left(\frac{\delta f\big(q^{-\frac 1 2}t\big)}{\delta t}
\right)^2d_qt\right)^{\frac 1 2}\,.
\end{array}
\end{equation}
Now, using the corresponding definitions of the $\,q$-integral and of the
operator $\,\delta\,$ one finds that
\begin{equation*}
\begin{array}{l}
\displaystyle\int_{-1}^{1}\left(\frac{\delta f\big(q^{\frac 1 2}t\big)}{\delta t}
\right)^2d_qt= \\ [1.2em]
\displaystyle (1-q)\sum_{n=0}^{\infty}\left\{\Big[f\big(q^n\big)-
f\big(q^{n+1}\big)\Big]^2+\Big[f\big(-q^n\big)-f\big(-q^{n+1}\big)\Big]^2
\right\}q^{-n}
\end{array}
\end{equation*}
hence, since $\,f\,$ is $\,q$-linear H\"older of order $\,\lambda>\frac 1 2\,,$
by (\ref{4.3}),
\begin{equation}\label{4.11}
\int_{-1}^{1}\left(\frac{\delta f\big(q^{\frac 1 2}t\big)}{\delta t}
\right)^2d_qt\leq 2M^2(1-q)\sum_{n=0}^{\infty}q^{(2\lambda-1)n}=
\frac{2(1-q)M^2}{1-q^{2\lambda-1}}\,.
\end{equation}
In a similar way we obtain
\begin{equation}\label{4.12}
\int_{-1}^{1}\left(\frac{\delta f\big(q^{-\frac 1 2}t\big)}{\delta
t} \right)^2d_qt\leq \frac{2(1-q)M^2}{1-q^{2\lambda-1}}\,.
\end{equation}
Thus, (\ref{4.9}) and (\ref{4.10}) become, respectively,
\begin{equation}\label{4.13}
\left|\int_{-1}^{1}S_q\left(q\omega_k t\right)
\frac{\delta f\big(q^{\frac 1 2}t\big)}{\delta t}d_qt\right|\leq \\ [1em]
\frac{\sqrt{2}q^{-\frac 1 4}(1-q)M}{\sqrt{1-q^{2\lambda-1}}}\left(C_q
\left(q^{\frac 1 2}\omega_k\right)S_q^{\prime}(\omega_k)\right)^{\frac 1 2}
\end{equation}
and
\begin{equation}\label{4.14}
\left|\int_{-1}^{1}C_q\left(q^{\frac 1 2}\omega_k t\right)
\frac{\delta f\big(q^{-\frac 1 2}t\big)}{\delta t}d_qt\right|\leq
\frac{\sqrt{2}(1-q)M}{\sqrt{1-q^{2\lambda-1}}}\left(C_q
\left(q^{\frac 1 2}\omega_k\right)S_q^{\prime}(\omega_k)\right)^{\frac 1 2}\,.
\end{equation}
Finally, using (\ref{4.13}) and (\ref{4.14}) in (\ref{4.4}) and (\ref{4.5}),
respectively, by Theorems A, E, F  and identity (\ref{2.10}),
as well as Lemma \ref{L3.1}, one concludes that the conditions of Theorem H are fulfilled with,
for instance, $\,c=3/2\,$, thus the $\,q$-Fourier series (\ref{3.1}), at
least for $\,0<q\leq (1/50)^{1/49}\,,$ converges uniformly on the set
$\,\displaystyle V_q=\left\{\pm q^{n-1}:\,n\in\bkN\,\right\},$  hence,
by Theorem I, under the same restriction on $\,q\,,$
$$\,S_q[f](x)=f(x)\;,\quad \forall\,x\in V_q=\left\{\pm q^{n-1}:\,n\in\bkN\,\right\}\,.$$
\end{proof}
A simple analysis of the previous theorem shows immediately that the behavior of
the function $\,f\,$ at the origin is crucial to study the convergence of the
$\,q$-Fourier series $\,S_q[f]\,.$ Consider, then, the following concept:

\vspace{0.7em}
\noindent$\mathbf{Definition\:4.2}$
\textsl{
A function $\,f\,$ is said to be almost $\,q$-linear H\"older of order $\,\lambda\,$
if two constants $\,M\,,$ $\,\lambda\,$ and a positive integer $\,n_0\,$
exist such that
\begin{equation}\label{4.15}
\Big|f\big(\pm q^{n-1}\big)-f\big(\pm q^{n}\big)\Big|\leq M q^{\lambda n}
\end{equation}
holds for every $\,n\geq n_0\,.$}

\vspace{0.4em}
Obviously that every {\em $\,q$-linear H\"older function of order $\,\lambda\,$}
is {\em almost $\,q$-linear H\"older function of order $\,\lambda\,$}.
\begin{cor}\label{C4.2}
If a function $\,f\in L_q^{\infty}[-1,1]\,$ is almost $\,q$-linear H\"older of
order $\,\lambda>\frac 1 2\,$ and satisfies $\,f(0^+)=f(0^-)\,$
then, at least for $\,0<q\leq (1/50)^{1/49}\,,$ the corresponding $\,q$-Fourier
series $\,S_q[f]\,$ converges uniformly to $\,f\,$ on the set of points
$\,\displaystyle V_q=\left\{\pm q^{n-1}:\,n\in\bkN\,\right\}\,.$
\end{cor}
\begin{proof} By hypothesis, $\,f\,$ is
almost $\,q$-linear H\"older of order $\,\lambda>1/2\,$,
i.e., it satisfies (\ref{4.15}). Then the relations (\ref{4.11}) and
(\ref{4.12})) now become
\begin{equation*}
\int_{-1}^{1}\left(\frac{\delta f\big(q^{\frac 1 2}t\big)}{\delta t}
\right)^2d_qt\leq \frac{2(1-q)M_1^2\,q^{n_0}}{1-q^{2\lambda-1}}
\end{equation*}
and
\begin{equation*}
\int_{-1}^{1}\left(\frac{\delta f\big(q^{-\frac 1 2}t\big)}{\delta
t} \right)^2d_qt\leq
\frac{2(1-q)M_2^2\,q^{n_0}}{1-q^{2\lambda-1}}\,,
\end{equation*}
respectively, where $\,M_1\,$ and $\,M_2\,$ are constants.
Therefore, using the above inequalities in formulas (\ref{4.9}) and
(\ref{4.10}) we get two new inequalities that differ from (\ref{4.13})
and (\ref{4.14}) only by a constant in the corresponding right hand
side. Hence, the conclusion on the uniform convergence follows.
\end{proof}

\begin{cor}\label{C4.3}
If $\,f\in L_q^{\infty}[-1,1]\,$ satisfies $\,f(0^+)=f(0^-)\,$ and there exists
a neighborhood of the origin where the function $\,f\,$ is continuous
and piecewise smooth then, at least for $\,0<q\leq (1/50)^{1/49}\,,$
the corresponding $\,q$-Fourier series $\,S_q[f]\,$ converges uniformly to
$\,f\,$ on the set of points
$\,\displaystyle V_q=\left\{\pm q^{n-1}:\,n\in\bkN\,\right\}\,.$
\end{cor}
\begin{proof} It's just a consequence of the fact that a function $\,f\,$ that is
continuous and piecewise smooth at any neighborhood  of the origin satisfies
a Lipschitz condition \cite[page 204]{N:1977}. Thus, it satisfies a H\"older
condition of order $1$ on that neighborhood and so, by Corollary \ref{C4.2}, the
uniform convergence follows.
\end{proof}

\section{Convergence on and outside the $q$-linear grid}

The convergence of the basic Fourier series (\ref{3.1})-(\ref{3.4})
always refer to the discrete set of the points of the $\,q$-linear grid
$\,\displaystyle \,V_q=\left\{\pm q^{n-1}:\,n\in\bkN\,\right\}$.

\noindent Two important questions arise at this moment:
\begin{itemize}
\item {\em The above mentioned $\,q$-Fourier series also converges
outside the points of the $\,q$-linear grid?}
\item {\em In that case, to what function it converges?}
\end{itemize}

Next theorem will give a positive answer to both questions.

\begin{teo}\label{T5.1}
Let $\,f\in L_q^{\infty}[-1,1]\,$ and suppose
that $\:c\in\bkR^+\:$ exists such that, as $\,k\rightarrow\infty\,$,
\begin{equation}
\int_{-1}^1\!f(t)C_q\!\left(q^{\frac 1 2}\omega_kt\right)\!d_qt=
\mathcal{O}\!\left(q^{(k+c)^2}\right)\,,\;\int_{-1}^1\!f(t)S_q\!\left(q\omega_kt\right)d_qt=
\mathcal{O}\!\left(q^{(k+c-\frac 1 2)^2}\right).
\label{5.1}
\end{equation}
If $\,f\,$ is analytic inside
$\,C_{\delta}=\left\{z\in\bkC\,:\:|z|<\delta\right\}\,,$ where
$\,\delta\,$ is a positive quantity such that $\,0<\delta\leq q^{-\sigma}\,$
with $\,0<\sigma<c\,,$ then, at least for
$\,0<q\leq \sqrt[50]{1/51}\,,$
\begin{equation}
f(z)=S_q[f](z)\quad\mbox{in}\quad C_{\delta}=\left\{\,z\in\bkC:\,|z|
<\,\delta\,\right\}\,.
\label{5.2}
\end{equation}
\end{teo}
\begin{proof} We first notice that
$$C_q\left(q^{\frac 1 2}\omega_kz\right)=\sum_{n=0}^{\infty}\frac{(-1)^n
q^{n(n-1)}}{\left(q^2,q;q^2\right)_n}q^{\frac 3 2 n}\omega_k^{2n}z^{2n}$$
and
$$S_q\left(q\omega_kz\right)=\frac{q\omega_kz}{1-q}\sum_{n=0}^{\infty}
\frac{(-1)^nq^{n(n-1)}}{\left(q^2,q^3;q^2\right)_n}q^{\frac 7 2 n}\omega_k^{2n}
z^{2n}$$
hence, for sufficiently large values of $\,k\,$, by Theorem A, whenever
$\,|z|\,\leq\,q^{-\sigma}\,,$
\begin{equation}\label{5.3}
\begin{array}{lll}
\left|C_q\left(q^{\frac 1 2}\omega_kz\right)\right| & \leq & \displaystyle\sum_{n=0}^{\infty}\frac{
q^{n(n-1)}}{\left(q^2,q;q^2\right)_n}q^{2n(1-k+\epsilon_k)}\left(q^{-\sigma}\right)^{2n} \\ [1em]
 & \leq & \displaystyle\frac{q^{-\left(k-\frac 1 2+\sigma-\epsilon_k\right)^2}}{(q;q)_{\infty}}
\sum_{n=0}^{\infty}q^{\left(n-k+\frac 1 2-\sigma+\epsilon_k\right)^2}
\end{array}
\end{equation}
and
\begin{equation}\label{5.4}
\begin{array}{lll}
\left|S_q\left(q\omega_kz\right)\right| & \leq & \displaystyle\frac{q\omega_kz}{1-q}
\sum_{n=0}^{\infty}\frac{q^{n(n-1)}}{\left(q^2,q^3;q^2\right)_n}q^{2n(2-k+\epsilon_k)}\left(q^{-\sigma}\right)^{2n} \\ [1em]
 & \leq & \displaystyle\frac{q^{\frac 5 4-k+\epsilon_k-\left(k-\frac 3 2+\sigma-\epsilon_k\right)^2}}{(q;q)_{\infty}}
\sum_{n=0}^{\infty}q^{\left(n-k+\frac 3 2-\sigma+\epsilon_k\right)^2}\,.
\end{array}
\end{equation}
An easy calculation shows that
\[
\begin{array}{lll}
\displaystyle \sum_{n=0}^{\infty}q^{\left(n-k+\frac 1 2+\epsilon_k-\sigma\right)^2}
& = & \displaystyle \sum_{n=0}^{k-1}q^{\left(n-k+\frac 1 2-\sigma+\epsilon_k\right)^2}+
\displaystyle \sum_{n=k}^{\infty}q^{\left(n-k+\frac 1 2-\sigma+\epsilon_k\right)^2} \\ [1em]
 & = & \displaystyle \sum_{m=0}^{k-1}q^{\left(m+\frac 1 2+\sigma-\epsilon_k\right)^2}+
\displaystyle \sum_{m=0}^{\infty}q^{\left(m+\frac 1 2-\sigma+\epsilon_k\right)^2}\,.
\end{array}
\]
thus, if
\begin{equation*} 
\left|\,\sigma\,\right|\,<\,\frac 1 2\,,
\end{equation*}
for sufficiently large values of $\,k\,,$
\begin{equation*} 
\begin{array}{lllllll}
\displaystyle \sum_{n=0}^{\infty}q^{\left(n-k+\frac 1 2+\epsilon_k-\sigma\right)^2}
 & < & \displaystyle \sum_{m=0}^{k-1}q^{m^2}+\sum_{m=0}^{\infty}q^{m^2}
 & < & \displaystyle 2\sum_{m=0}^{\infty}q^{m} & = & \displaystyle\frac{2}{1-q}\,.
\end{array}
\end{equation*}
In a similar way, for a given $\;p\in\bkN_0\;$, if
\begin{equation}
\left|\,\sigma\,\right|\,<\,\frac 1 2+p
\label{5.7}
\end{equation}
then, for sufficiently large values of $\,k\,,$
\begin{equation}
\begin{array}{lll}
\displaystyle \sum_{n=0}^{\infty}q^{\left(n-k+\frac 1 2+\epsilon_k-\sigma\right)^2}
 & < & \displaystyle 2p+\frac{2}{1-q}\,.
\end{array}
\label{5.8}
\end{equation}
With the same reasoning we get, again for sufficiently large values of $\,k\,,$
\begin{equation}
\begin{array}{lll}
\displaystyle \sum_{n=0}^{\infty}q^{\left(n-k+\frac 3 2+\epsilon_k-\sigma\right)^2}
 & < & \displaystyle 2p+\frac{2}{1-q}\,.
\end{array}
\label{5.9}
\end{equation}
Hence, by (\ref{5.3}), (\ref{5.8}) and (\ref{5.4}), (\ref{5.9}), we may write,
respectively, for $\,k\,$ large enough,
\begin{equation}\label{5.10}
\left|C_q\left(q^{\frac 1 2}\omega_kz\right)\right| \leq \frac{2p(1-q)+2}{(q;q)_{\infty}}
\,q^{-\left(k-\frac 1 2+\sigma-\epsilon_k\right)^2}
\end{equation}
and
\begin{equation}\label{5.11}
\left|S_q\left(q\omega_kz\right)\right| \leq \frac{2p(1-q)+2}{(q;q)_{\infty}}
\,q^{\frac 5 4-k+\epsilon_k-\left(k-\frac 3 2+\sigma-\epsilon_k\right)^2}\,.
\end{equation}
This way, for $\,k\,$ large enough, using (\ref{3.2}) and (\ref{3.4}),
Theorems E and F, relation (\ref{2.10}) and inequality (\ref{5.10}), at least
for $\,0<q\leq \sqrt[50]{1/51}\,,$
\begin{equation*}
\left|a_kC_q\!\left(q^{\frac 1 2}\omega_kz\right)\right|\leq
\frac{2p(1\!-\!q)\!+\!2}{(1\!-\!q)^2(q;q)_{\infty}^2}
\left|\int_{-1}^1\!f(t)C_q\!\left(q^{\frac 1 2}\omega_kt\right)d_qt\right|
\frac{q^{-\left(k-\frac 1 2+\sigma-
\epsilon_k\right)^2-k+\frac 1 4+\epsilon_k}}{\left|S_k\right|}.
\end{equation*}
By hypothesis (\ref{5.1}), we may suppose that $\,c_1\in\bkR^+\,$ and $\,M_1>0\,$
exist such that, for $\,k\,$ large enough,
\begin{equation}
\left|\int_{-1}^1f(t)C_q\left(q^{\frac 1 2}\omega_kt\right)d_qt\right|\leq
M_1q^{(k+c_1)^2}\,.
\label{5.12}
\end{equation}
In that case we have
\begin{equation*} 
\!\left|a_kC_q\left(q^{\frac 1 2}\omega_kz\right)\right|\leq 2M_1
\frac{p(1\!-\!q)+1}{(\!1\!-\!q)^2(q;q)_{\infty}^2}
\frac{q^{(k+\frac{c_1+\sigma}{2}-\frac 1 4-\frac{\epsilon_k}{2})(1+2(c_1-\sigma)+
2\epsilon_k)-k+\frac 1 4+\epsilon_k}}{\left|S_k\right|}
\end{equation*}
hence, if $\,1+2(c_1-\sigma)\,>\,1\,,$ i.e., if $\,\sigma\,<\,c_1\,$
then, taking into account Theorem A and (\ref{2.11}), and the Theorems E
and F, at least for $\,0<q\leq \sqrt[50]{1/51}\,,$
\begin{equation}
\left|a_kC_q\left(q^{\frac 1 2}\omega_kz\right)\right|\leq A_1q^{\theta_1k}\quad,
\label{5.14}
\end{equation}
where $\,A_1\,$ and $\,\theta_1\,$ are positive constants.

Analogously, for $\,k\,$ large enough, (\ref{3.3}) and (\ref{3.4}),
Theorems E and F, relation (\ref{2.10}) and inequality (\ref{5.11}),
\begin{equation*}
\left|b_kS_q\left(q\omega_kz\right)\right|\leq
\frac{2p(1-q)+2}{(1-q)^2(q;q)_{\infty}^2}
\left|\int_{-1}^1f(t)S_q\left(q\omega_kt\right)d_qt\right|
\frac{q^{-\left(k-\frac 3 2+\sigma-\epsilon_k\right)^2-2k+2+2\epsilon_k}}{
\left|S_k\right|}
\end{equation*}
so, again by hypothesis (\ref{5.1}), if we admit that $\,c_2\in\bkR^+\,$ and
$\,M_2>0\,$ exist such that
\begin{equation}
\left|\int_{-1}^1f(t)S_q\left(q\omega_kt\right)d_qt\right|\leq
M_2q^{(k+c_2-\frac 1 2)^2}\,,
\label{5.15}
\end{equation}
then,
\begin{equation*} 
\hspace{-0.3em}\left|b_kS_q\left(q\omega_kz\right)\right|\leq 2M_2
\frac{p(1\!-\!q)+1}{(1\!-\!q)^2(q;q)_{\infty}^2}
\frac{q^{(k+\frac{c_2+\sigma}{2}-\frac 3 4-
\frac{\epsilon_k}{2})(2+2(c_2-\sigma)+2\epsilon_k)-2k+2+2\epsilon_k}}{
\left|S_k\right|}\,.
\end{equation*}
Similarly, if $\,2+2(c_2-\sigma)\,>\,2\,,$ i.e., if $\,\sigma\,<\,c_2\,$
then, at least for $\,q\,$ such that $\,0<q\leq \sqrt[50]{1/51}\,,$
\begin{equation}
\left|b_ks_q\left(q\omega_kz\right)\right|\leq A_2q^{\theta_2k}\,,
\label{5.17}
\end{equation}
being $\,A_2\,$ and $\,\theta_2\,$ positive constants.

We remark that in (\ref{5.7}) we may choose $\,p\,$ sufficiently large
in order that one haves
\begin{equation}
\;-\frac 1 2-p\,<\,0\,<\,\sigma\,<\,\min\left\{c_1,c_2\right\}\,\leq\,\frac 1 2+p\,,
\label{5.18}
\end{equation}
thus, replacing $\,c_1\,$ and $\,c_2\,$ from (\ref{5.12}) and
(\ref{5.15}) by $\,c=\min\left\{c_1,c_2\right\}\,$, respectively,
we conclude, through (\ref{5.14}) and (\ref{5.17}), that the conditions (\ref{5.1})
guaranty the uniform convergence of the $\,q$-Fourier series (\ref{3.1}) in
$\,C_{q^{-\sigma}}=\left\{\,z\in\bkC\,:\,|z|\,<\,q^{-\sigma}\,\right\}\,$
if $\,\sigma\,$ satisfies (\ref{5.18}).
This way, under this condition on $\,\sigma\,$, we have, by Theorem H,
$$f(x)=S_q[f](x)\quad\mbox{whenever}\quad x\in V_q\,,$$
since $\,V_q\subset C_{q^{-\sigma}}\,,$ where
$\,V_q=\left\{\,q^{n-1}:\,n\in\bkN\right\}\,$
is the corresponding set of Theorem I and $\,C_{q^{-\sigma}}\,$
is the interior of the circle of the complex plane with center at the origin
and radius $\,q^{-\sigma}\,.$

On the other side, again by the uniform convergence of the $\,q$-Fourier series
$\,S_q[f](x)\,$ on $\,C_{q^{-\sigma}}\,,$ since the terms of the mentioned
$\,q$-Fourier series are entire functions we then have that the
$\,q$-series is analytic inside $\,C_{q^{-\sigma}}\,.$
From the continuity of both members of the above equality it results
$\,f(0)=S_q[f](0)\,.$ Thus, if $\,f\,$ is analytic inside
$\,C_{\delta}=\left\{\,z\in\bkC\,:\,|z|<\delta\,\right\}\,,$ where
$\,0\,<\,\delta\,\leq\,q^{-\sigma}\,,$ then $\,f(z)\,$ and $\,S_q[f](z)\,$
are analytic inside $\,C_{\delta}\,$ and coincide in a set with a limit point
in the interior of such circle; by the {\em principle of analytic continuation}
\cite[Corollary 4.4.1]{D:1984}, the above mentioned functions must coincide
in the whole set $\,C_{\delta}\,,$ which proves (\ref{5.2}).

\end{proof}

\section{Examples}
In this section we will present four examples of $q$-Fourier
series and study the corresponding questions about convergence.

\vspace{1em}
{\em Example 1}: $g(x)=|x|$

\vspace{0.6em}
\noindent The basic Fourier series of the absolute value function is
given \cite{C:2005} by
\begin{equation*}
S_q[g](x)=\displaystyle\frac{1}{1+q}-2q^{-\frac 1
2}(1-q)\sum_{k=1}^{\infty}\frac{1-C_q\left(q^{\frac 1
2}\omega_k\right)} {\omega_k^2C_q\left(q^{\frac 1
2}\omega_k\right)S_q^{\prime}(\omega_k)}C_q\left(q^{\frac 1
2}\omega_kx\right)\:.
\end{equation*}
Conditions of Theorem H are fulfilled \cite{C:2005} with, for
instance, $\,c=2\,.$ Thus, at least for $\,0<q\leq (1/50)^{1/49}\,,$
the $\,q$-Fourier series of the function $\,f(x)=\left|x\right|\,$ converges
uniformly on the set $\,V_q=\left\{\,\pm q^{n-1}:\,n\,\in\bkN\,\right\}\,$ so,
under the same restrictions on $\,q\,,$ by Theorem I,
\begin{equation*}
|x|\,=\,\frac{1}{1+q}-2q^{-\frac 1
2}(1-q)\sum_{k=1}^{\infty}\frac{1-C_q\left(q^{\frac 1
2}\omega_k\right)} {\omega_k^2C_q\left(q^{\frac 1
2}\omega_k\right)S_q^{\prime}(\omega_k)}C_q\left(q^{\frac 1
2}\omega_kx\right)\:
\end{equation*}
for all $\,x\in V_q=\left\{\,\pm q^{n-1}:\,n\,\in\bkN\,\right\}\,.$

\vspace{0.5em}
Now, we may obtain the same conclusion in a easier way through Theorem
\ref{T4.1}, by simple arguing that the absolute value function
\begin{itemize}
\item is bounded on $\,V_q=\left\{\,\pm q^{n-1}:\,n\,\in\bkN\,\right\},$
\item is continuous at the origin,
\item and satisfies the $\,q$-linear H\"older condition of
order $\,1\,$ since
$$\Big| \big|\pm q^{n-1}\big|-\big|\pm q^{n}\big|\Big|\leq (1-q)q^{n-1}\,.$$
\end{itemize}
Thus, by Theorem \ref{T4.1}, the same conclusion over the uniform convergence
follows. Notice that Corollaries \ref{C4.2} or \ref{C4.3} also apply.

\vspace{1.2em}
Given a function $\,f\,$, it is important to point out that Theorem \ref{T4.1}
or its Corollaries \ref{C4.2} and \ref{C4.3}, enable one to decide over the
uniform convergence of the $\,q$-Fourier series $\,S_q[f]\,$ without the need
to compute the corresponding coefficients: only requires a short study of the
function itself.

\vspace{1.8em}
{\em Example 2}: \hspace{0.5em} $h(x)=\left\{
\begin{array}{lll}
-1 & \mbox{if} & x\leq 0 \\
  & & \\
1 & \mbox{if} & x>0
\end{array}\right.$

\vspace{0.8em}
\noindent In this example, the conditions of Theorem H were not satisfied
\cite[Remark 3]{C:2005}. It was shown, using Theorem G, that the
$\,q$-Fourier series
\begin{equation*}
S_q[h](x)=\displaystyle
2\sum_{k=1}^{\infty}\frac{1-C_q\left(q^{\frac 1 2}\omega_k\right)}
{\omega_kC_q\left(q^{\frac 1
2}\omega_k\right)S_q^{\prime}(\omega_k)}S_q(q\omega_kx)
\end{equation*}
is (pointwise) convergent at each (fixed) point $\,x\in V_q\,$.
Theorem \ref{T4.1} doesn't apply too (neither its corollaries) since
$\,h\big(0^+\big)\neq h\big(0^-\big)\,.$

\vspace{1em} {\em Example 3}: $\:H^{(a)}(x)=\left\{
\begin{array}{lll}
-1 & \mbox{se} & x\leq a \\
  & & \hspace{4em};\hspace{2em}(a>0)\\
1 & \mbox{se} & x>a
\end{array}\right.
\:$
\vspace{0.6em}
\newline
Once $\,0<q<1\,$ is fixed, denote by $\,n_a\,$ the least positive integer $\,j\,$
such that $\,q^j<a\,,$ i.e., $\,n_a=\left[\log_q a\right]+1\,.$ Then
\begin{equation}\label{6.1}
a_0 = -2q^{n_a}
\end{equation}
and, for $\:k=1,2,3,\ldots\,,$
\begin{equation*}
a_k=\frac{2(1-q)}{q^{-\frac 1 2+n_a}\omega_k^2\mu_k}\left[C_q\left(q^{\frac 1 2+n_a}\omega_k\right)-
C_q\left(q^{-\frac 1 2+n_a}\omega_k\right)\right]\,.
\end{equation*}
By Theorem D, $$C_q\left(q^{\frac 1 2+n_a}\omega_k\right)-
C_q\left(q^{-\frac 1 2+n_a}\omega_k\right)=
q^{-\frac 1 2+n_a}\omega_kS_q\left(q^{n_a}\omega_k\right)\,,$$
thus
\begin{equation}
a_k=-\frac{2(1-q)S_q\left(q^{n_a}\omega_k\right)}{\omega_k\mu_k}=
-\frac{2}{\omega_k}\frac{S_q\left(q^{n_a}\omega_k\right)}
{C_q\left(q^{\frac 1 2}\omega_k\right)S_q^{\prime}(\omega_k)}\,.
\label{6.2}
\end{equation}
For $\,k=1,2,3,\ldots\,$ we have
\begin{equation*}
b_k = -\displaystyle\frac{2(1-q)}{\omega_k^2\mu_k}\left[\frac{S_q\left(q^{1+n_a}\omega_k\right)-
S_q\left(q^{n_a}\omega_k\right)}{q^{n_a}}-S_q(q\omega_k)\right]\,.
\end{equation*}
By Theorem D,
$$S_q\left(q^{1+n_a}\omega_k\right)-S_q\left(q^{n_a}\omega_k\right)=
-q^{n_a}\,\omega_k\,C_q\!\left(q^{\frac 1 2+n_a}\omega_k\right)\,,$$
so, by (\ref{2.10}),
\begin{equation}
b_k=\displaystyle\frac{2(1\!-\!q)}{\omega_k\mu_k}
\left[C_q\!\left(q^{\frac 1 2+n_a}\omega_k\!\right)\!-C_q\!\left(q^{\frac 1 2}\omega_k\!\right)\right]=
\frac{2}{\omega_k}\frac{C_q\!\left(q^{\frac 1 2+n_a}\omega_k\!\right)\!-C_q\!\left(q^{\frac 1 2}\omega_k\!\right)}
{C_q\!\left(q^{\frac 1 2}\omega_k\right)S_q^{\prime}(\omega_k)}\,.
\label{6.3}
\end{equation}
hence, substituting  (\ref{6.1}), (\ref{6.2}) and (\ref{6.3}) into
(\ref{3.1}) it becomes
\begin{equation}
\begin{array}{l}
S_q[H^{(a)}](x)=-q^{n_a}- \\ [1em]
\displaystyle\hspace{1.5em}2\sum_{k=1}^{\infty}\frac{S_q\left(q^{n_a}\omega_k\right)C_q\!\left(q^{\frac 1 2}\omega_kx\right)+
\left[C_q\!\left(q^{\frac 1 2}\omega_k\right)-C_q\!\left(q^{\frac 1 2+n_a}\omega_k\right)\right]S_q(q\omega_kx)}
{\omega_k\,C_q\!\left(q^{\frac 1 2}\omega_k\right)S_q^{\prime}(\omega_k)}\,.
\end{array}
\label{6.4}
\end{equation}
We notice that {\em Example 2}\: follows from {\em Example 4}\: by computing the limit
$\,n_a\rightarrow\infty\,,$ i.e., when $\,a\rightarrow 0\,.$
Again by Theorem D,
$$S_q(q^{n_a}\omega_k)=S_q(q\omega_k)\sum_{j=0}^{n_a-1}(-1)^{j}q^{j(j+\frac{1}{2})}
\frac{\left(q^{n_a-j};q\right)_{2j+1}}{(q;q)_{2j+1}}\omega_k^{2j}$$
and
$$C_q(q^{\frac 1 2+n_a}\omega_k)=C_q(q^{\frac 1 2}\omega_k)\sum_{j=0}^{n_a}
(-1)^{j}q^{j(j-\frac{1}{2})}\frac{\left(q^{1+n_a-j};q\right)_{2j}}{(q;q)_{2j}}\omega_k^{2j}\,,$$
thus, since $\,S_q(q\omega_k)=-\omega_kC_q\big(q^{1/2}\omega_k\big)\,,$
for $\,k=1,2,3,\ldots\,,$
\begin{equation*}
\int_{-1}^1\!H^{(a)}(x)C_q(q^{\frac 1 2}\omega_k x)d_qt=2(1\!-\!q)
C_q\!\left(q^{\frac{1}{2}}\omega_k\!\right)\!\sum_{j=0}^{n_a-1}\!(-1)^{j}q^{j(j+\frac{1}{2})}
\frac{\left(q^{n_a-j};q\right)_{2j+1}}{(q;q)_{2j+1}}\omega_k^{2j}
\end{equation*}
and
\begin{equation*}
\begin{array}{l}
\displaystyle\int_{-1}^{1}H^{(a)}(x)S_q(q\omega_k x)d_qt=2q^{-\frac
1 2}(1-q)
\frac{c_q\left(q^{\frac{1}{2}}\omega_k\right)}{\omega_k}\times \\
[1em]
\hspace{10em}\displaystyle\left[\:\sum_{j=0}^{n_a}(-1)^{j}q^{j(j-\frac{1}{2})}
\frac{\left(q^{1+n_a-j};q\right)_{2j}}{(q;q)_{2j}}\omega_k^{2j}-1\:\right]\,.
\end{array}
\end{equation*}
For each fixed $\,a>0\,$, at least for $\,0<q\leq (1/50)^{1/49}\,,$
the $\,q$-Fourier series (\ref{6.4}) converges uniformly on the set
$\,V_q=\left\{\,\pm q^{n-1}:\,n\,\in\bkN\,\right\}\,$: in fact, after some
computations, one verifies that the conditions of Theorem H are satisfied
with, for instance, $\,c=2\,,$ hence, whenever $\,x\in V_q\,$ and under the
above restriction on $\,q\,,$ we may write by Theorem I,
\begin{equation}\label{6.5}
\begin{array}{l}
\displaystyle H^{(a)}(x)\equiv -q^{n_a}- \\ [1em]
\hspace{1em}\displaystyle 2\sum_{k=1}^{\infty}\frac{S_q\left(q^{n_a}\omega_k\right)C_q\left(q^{\frac 1 2}\omega_kx\right)+
\left[C_q\left(q^{\frac 1 2}\omega_k\right)-C_q\left(q^{\frac 1 2+n_a}\omega_k\right)\right]S_q(q\omega_kx)}
{\omega_kC_q\left(q^{\frac 1 2}\omega_k\right)S_q^{\prime}(\omega_k)}\,.
\end{array}
\end{equation}

Another approach is the following: one easily check that
$\,H^{(a)}\in L_q^{\infty}[-1,1]\,,$ $\,H^{(a)}\big(0^+\big)=0=H^{(a)}\big(0^-\big)\,$
and $\,H^{(a)}\,$ is almost $\,q$-linear H\"older of order bigger then $\,\frac 1 2\,$
since $$\left|H^{(a)}\big(\pm q^{n-1}\big)-H^{(a)}\big(\pm q^{n}\big)\right|=0
\;,\quad n\geq n_a+1=\left[\log_q a\right]+2\,.$$
By Corollary \ref{C4.2}, the $\,q$-Fourier series $\,S_q\big[H^{(a)}\big]\,$
converges uniformly on the set $\,V_q\,,$ thus (\ref{6.5}) follows.

\vspace{1em} {\em Example 4}: $f(x)=x^m$

\vspace{0.5em}
\noindent In \cite[Proposition 6.1]{C:2005} it was presented the Fourier
expansion of the function $\:f(x)=x^m\,,$ $m=0,1,2,\ldots\,,$ in terms
of the functions $\,C_q\,$ and $\,S_q\,:$
\begin{equation*}
\begin{array}{l}
S_q[x^m](x)\,=\,\displaystyle\frac{1+(-1)^m}{2}\frac{1-q}{1-q^{m+1}}+ \\ [0.8em]
\hspace{3em}(q;q)_m{\displaystyle\sum_{k=1}^{\infty}}\left\{\frac{1\!+\!(-1)^m}{S_q^{\prime}(\omega_k)}
{\displaystyle\sum_{i=0}^{\left[\frac{m-2}{2}\right]}}\displaystyle\frac{(-1)^iq^{(i+1)(i-m+\frac 1 2)}}{
\omega_k^{2i+2}(q;q)_{m-1-2i}}C_q(q^{\frac 1 2}\omega_kx)\;+\right. \\ [1em]
\hspace{3em}\displaystyle \left.q^{\frac 1 2}\frac{(-1)\!+\!(-1)^m}{S_q^{\prime}(\omega_k)}
{\displaystyle\sum_{i=0}^{\left[\frac{m-1}{2}\right]}}
\frac{(-1)^iq^{(i+1)(i-m-\frac 1 2)}}{\omega_k^{2i+1}(q;q)_{m-2i}}S_q(q\omega_kx)\right\}\,,
\end{array}
\end{equation*}
where $\:[x]\:$ denotes the greatest integer which does not
exceed $\:x\:$ and we will take as zero a sum where the superior
index is less than the inferior one.

\noindent Furthermore, it was proved that the conditions of Theorem H are
fulfilled with , for instance, $\,c=2\,.$ Thus, at least for $\,0<q\leq (1/50)^{1/49}\,,$ the
$\,q$-Fourier series of the function $\,f(x)=x^m\,$ converges
uniformly on the set $\,V_q=\left\{\,\pm q^{n-1}:\,n\,\in\bkN\,\right\}\,,$ so,
by Theorem I,
$$ x^m\,=\,S_q[x^m](x)\quad \mbox{whenever}
\quad x\in V_q=\left\{\,\pm q^{n-1}:\,n\,\in\bkN\,\right\}\,.$$
We notice that the conditions of Theorem \ref{T4.1} are trivial
checked when $\,f(x)=x^m\,.$

\vspace{0.2em}
Now, since $\,f\,$ satisfies the conditions of Theorem \ref{T5.1} with, for instance,
$\,c=1\,$ and $\,f\,$ is an entire function then, by Theorem \ref{T5.1},
$$S_q[x^m](x)=x^m\:,\quad \forall x\in C_{\delta}=\left\{\,z\in\bkC:\,|z|
<\,\delta\,\right\}$$
where $\,0<\delta<q^{-\sigma}\,$ and $\,0<\sigma<1\,.$


\subsection*{Concluding remarks}
We notice that Theorem \ref{T4.1} or Corollaries \ref{C4.2} and \ref{C4.3} are
$\,q\,$-analogs of the corresponding classical theorems on uniform convergence for
trigonometric Fourier series. See, for instance, Theorem 1 of \cite[page 204]{N:1977}
or Theorem 55 of \cite[page 41]{HR:1999}.

Mathematica$^{\mbox{\scriptsize\copyright}}$ suggests that Theorems
(\ref{4.1}) and (\ref{5.1}) remain valid for $\,0<q<1\,$.
It's an open question and to prove it a different technic is required.

\subsection*{Acknowledgements}
Discussions with Prof. R. \'Alvarez-Nodarse
are kindly acknowledged. This research was partially supported by CMUC
from the University of Coimbra.


\end{document}